\documentclass[12pt]{amsart}
\usepackage{amssymb,amsmath,amsfonts,latexsym,setspace}
\usepackage{bm}
\newcommand{\newsection}[1]{\setcounter{equation}{0} \section{#1}}
\newcommand{\bea}{\begin{eqnarray}}
\newcommand{\eea}{\end{eqnarray}}

\newcommand{\vp}{\varphi}

\newcommand{\cld}{\mathcal{D}}
\newcommand{\cle}{\mathcal{E}}

\newcommand{\clh}{\mathcal{H}}
\newcommand{\clk}{\mathcal{K}}
\newcommand{\cll}{\mathcal{L}}

\newcommand{\raro}{\rightarrow}

\def \qed {\hfill \vrule height6pt width 6pt depth 0pt}
\def\textmatrix#1&#2\\#3&#4\\{\bigl({#1 \atop #3}\ {#2 \atop #4}\bigr)}
\def\dispmatrix#1&#2\\#3&#4\\{\left({#1 \atop #3}\ {#2 \atop #4}\right)}
\newcommand{\be}{\begin{equation}}
\newcommand{\ee}{\end{equation}}
\newcommand{\ben}{\begin{eqnarray*}}
\newcommand{\een}{\end{eqnarray*}}

\newcommand{\NI}{\noindent}

\newcommand{\bi}{\begin{itemize}}
\newcommand{\ei}{\end{itemize}}

\newtheorem{Theorem}{\sc Theorem}[section]
\newtheorem{Lemma}[Theorem]{\sc Lemma}
\newtheorem{Proposition}[Theorem]{\sc Proposition}
\newtheorem{Corollary}[Theorem]{\sc Corollary}
\newtheorem{Definition}[Theorem]{\sc Definition}
\newtheorem{Example}[Theorem]{\sc Example}
\newtheorem{Remark}[Theorem]{\sc Remark}
\newtheorem{Note}[Theorem]{\sc Note}
\newtheorem{Question}{\sc Question}
\newtheorem{ass}[Theorem]{\sc Assumption}
\newcommand{\bt}{\begin{Theorem}}
\def\beginlem{\begin{Lemma}}
\def\beginprop{\begin{Proposition}}
\def\begincor{\begin{Corollary}}
\def\begindef{\begin{Definition}}
\def\beginexamp{\begin{Example}}
\def\beginrem{\begin{Remark}}
\def\beginq{\begin{Question}}
\def\beginass{\begin{ass}}
\def\beginnote{\begin{Note}}
\newcommand{\et}{\end{Theorem}}
\def\endlem{\end{Lemma}}
\def\endprop{\end{Proposition}}
\def\endcor{\end{Corollary}}
\def\enddef{\end{Definition}}
\def\endexamp{\end{Example}}
\def\endrem{\end{Remark}}
\def\endq{\end{Question}}
\def\endass{\end{ass}}
\def\endnote{\end{Note}}

\begin{document}

\title[Curvature invariant and generalized
canonical Operator models - I]{Curvature invariant and generalized
canonical Operator models - I}
\author[Douglas]{Ronald G. Douglas}
\address{Department of Mathematics\\ Texas A\&M University\\ College Station\\ TX 77843\\ USA}
\email{rdouglas@math.tamu.edu}

\author[Kim]{Yun-Su Kim}
\address{Department of Mathematics\\The University of Toledo\\Toledo, OH 43606\\USA\\\textsf{Deceased}}

\author[Kwon]{Hyun-Kyoung Kwon}
\address{Department of Mathematical Sciences\\ Seoul National University\\ Seoul, 151-747\\ Republic of Korea}
\email{hyunkwon@snu.ac.kr}

\author[Sarkar]{Jaydeb Sarkar}
\address{Indian Statistical Institute\\
Statistic and Mathematics Unit\\ 8th Mile, Mysore Road\\ Bangalore 560059\\
India}
\email{jay@isibang.ac.in}

\thanks{The work of Douglas and Sarkar was partially supported by a grant from the National Science
Foundation. The work of Kwon was supported by the Basic Science
Research Program through the National Research Foundation of Korea
(NRF) grant funded by the Korean government (MEST) (No.
2010-0024371), and in part, by a Young Investigator Award at the
NSF-sponsored Workshop in Analysis and Probability, Texas A \& M
University, 2009. The research began in the summer of 2009. Sarkar
was at Texas A \& M University at the time and Kim and Kwon were
participants of the workshop. Sarkar would also like to acknowledge the hospitality of the mathematics departments of Texas A \& M University and the University of Texas at San Antonio, where part of his research was done.}



\begin{abstract}
One can view contraction operators given by a canonical model of Sz.-Nagy and Foias
as being defined by a quotient module where the basic building blocks are Hardy spaces. In
this note we generalize this framework to allow the Bergman and weighted Bergman spaces as
building blocks, but restricting attention to the case in which the operator obtained is in the
Cowen-Douglas class and requiring the multiplicity to be one. We view the classification of such
operators in the context of complex geometry and obtain a complete classification up to unitary
equivalence of them in terms of their associated vector bundles and their curvatures.

\end{abstract}

\maketitle

\newsection{Introduction}
One goal of operator theory is to obtain unitary invariants,
ideally, in the context of a concrete model for the operators
being studied. For a multiplication operator on a space of
holomorphic functions on the unit disk $\mathbb{D}$, which happens
to be contractive, there are two distinct approaches to models and
their associated invariants, one due to Sz.-Nagy and Foias
\cite{NF} and the other due to M. Cowen and the first author
\cite{CD}. The starting point for this work was an attempt to
compare the two sets of invariants and models obtained in these
approaches. We will work at the simplest level of generality
for which these questions make sense. Extensions of these results to more
general situations are pursued later in \cite{DKKS}.

For the Sz.-Nagy-Foias canonical model theory, the Hardy space
$H^2=H^2(\mathbb{D})$, of holomorphic functions on the unit disk
$\mathbb{D}$ is central if one allows the functions to take values
in some separable Hilbert space $\cle$. In this case, we will
now denote the space by $H^2 \otimes \cle$. One can view the
canonical model Hilbert space (in the case of a $C_{\cdot 0}$
contraction $T$) as given by the quotient of $H^2 \otimes \cle_*$,
for some Hilbert space $\cle_*$, by the range of a map
$M_{\Theta}$ defined to be multiplication by a contractive
holomorphic function, $\Theta(z) \in \cll(\cle, \cle_*)$, from
$H^2 \otimes \cle$ to $H^2 \otimes \cle_*$. If one assumes that
the multiplication operator associated with $\Theta(z)$ defines an
isometry (or is inner) and $\Theta(z)$ is purely contractive, that is, $\|\Theta(0) \eta\| < \|\eta\|$ for all $\eta (\neq 0)$ in $\cle$, then $\Theta(z)$ is the characteristic
operator function for the operator $T$. Hence, $\Theta(z)$
provides a complete unitary invariant for the compression of
multiplication by $z$ to the quotient Hilbert space of $H^2
\otimes \cle_*$ by the range of $\Theta(z)$. In general, neither
the operator $T$ nor its adjoint $T^*$ is in the $B_n(\mathbb{D})$
class of \cite{CD} but we are interested in the case in which the
adjoint $T^*$ is in $B_n(\mathbb{D})$ and we study the relation
between its complex geometric invariants (see \cite{CD}) and
$\Theta(z)$.

We use the language of Hilbert modules \cite{DP} which we believe to be
natural in this context. The Cowen-Douglas theory can also be recast in
the language of Hilbert modules \cite{CD1}. With this approach, the problem of the unitary equivalence of operators becomes identical
to that of the isomorphism of the corresponding Hilbert modules.

Furthermore, we consider ``models" obtained as quotient Hilbert modules in
which the Hardy module is replaced by other Hilbert modules of
holomorphic functions on $\mathbb{D}$ such as the Bergman module $A^2=A^2(\mathbb{D})$ or
the weighted Bergman modules $A_{\alpha}^2=A_{\alpha}^2(\mathbb{D})$ with weight parameter $\alpha > -1$.
We require in these cases that some analogue of the corona condition holds for the multiplier $\Theta(z)$.

As previously mentioned, we concentrate on a particularly simple case of the problem. We focus on the case of $\Theta \in H ^{\infty} _{\cll(\mathbb{C}, \mathbb{C}^2)}$, where $H ^{\infty} _{\cll(\mathbb{C}, \mathbb{C}^2)} = H ^{\infty} _{\cll(\mathbb{C}, \mathbb{C}^2)}(\mathbb{D})$ is the space of bounded, holomorphic $\cll(\mathbb{C}, \mathbb{C}^2)$-valued functions on $\mathbb{D}$, so that $\Theta(z)=\theta_1(z)\otimes e_1+
\theta_2(z)\otimes e_2$ for an orthonormal basis $\{e_1,e_2\}$ for $\mathbb{C}^2$ and $\theta_i(z) \in \cll(\mathbb{C})$, $i=1,2,$ and $z \in \mathbb{D}$. We shall adopt the notation
$\Theta=\{\theta_1, \theta_2 \}$. Recall that $\Theta$ is
said to satisfy the {\it corona condition} if there exists an
$\epsilon > 0$ such that $$|\theta_1(z)|^2 +
|\theta_2(z)|^2 \geq \epsilon,$$ for all $z \in \mathbb{D}$.
Moreover, we will use the notation $\clh_\Theta$ to denote the
quotient Hilbert module $(\clh \otimes \mathbb{C}^2) / \Theta
\clh$, where $\clh$ is the Hardy, the Bergman, or a weighted Bergman module.

Now we state the main results in this note which we will prove in Section 4. Let $\Theta=\{\theta_1, \theta_2\}$ and $\Phi=\{\vp_1, \vp_2\}$ both satisfy
the corona condition and denote by $\bigtriangledown^2$
the Laplacian $\bigtriangledown^2=4\partial
\bar{\partial}=4\bar{\partial}{\partial}$.\\

\NI \textsf{Theorem 4.4. }
The quotient Hilbert
modules $\clh_{\Theta}$ and $\clh_{\Phi}$ are isomorphic if and only
if
$$
\bigtriangledown^2 \mbox{log}\,\frac{|\theta_1(z)|^2 +
|\theta_2(z)|^2}{|\vp_1(z)|^2 + |\vp_2(z)|^2}=0,
$$
for all $z \in \mathbb{D}$, where $\clh$ is the Hardy module $H^2$,
the Bergman module $A^2$, or a weighted Bergman module $A^2_{\alpha}$. \\

\NI \textsf{Theorem 4.5. }
The quotient Hilbert modules
$(A^2_{\alpha})_{\Theta}$ and
$(A^2_{\beta})_{\Phi}$ are isomorphic if and
only if $\alpha=\beta$ and
$$
\bigtriangledown^2 \mbox{log}\,\frac{|\theta_1(z)|^2 +
|\theta_2(z)|^2}{|\vp_1(z)|^2 + |\vp_2(z)|^2}=0,
$$
for all $z \in \mathbb{D}$. \\

\NI \textsf{Theorem 4.7. }
Under no circumstances can $(H^2)_{\Theta}$
be isomorphic to $(A^2_{\alpha})_{\Phi}$.\\

\newsection{Hilbert Modules}

In the present section and the next, we take care of some preliminaries. We
begin with the following definition.

\begin{Definition}
Let $T$ be a linear operator on a Hilbert space $\clh$. We say that $\clh$ is a contractive Hilbert module over $\mathbb{C}[z]$ relative to $T$ if the module action from $\mathbb{C}[z] \times \clh$ to $\clh$ given by $$p \cdot f \mapsto p(T) f,$$ for $p \in \mathbb{C}[z]$ defines bounded operators such that $$\|p \cdot f\|_{\clh} = \|p(T) f\|_{\clh} \leq \|p\|_{\infty} \|f\|_{\clh},$$ for all $f \in \clh$, where $\|p\|_{\infty}$ is the supremum norm of $p$ on $\mathbb{D}$.
\end{Definition}

The module multiplication by the coordinate function will be denoted by $M_z$, that is, \[M_z f = z \cdot f = T f,\]for all $f \in \clh$.
\begin{Definition}
Given two Hilbert modules $\clh$ and $\tilde{\clh}$ over $\mathbb{C}[z]$, we say that $X: \clh \raro \tilde{\clh}$ is a module map if it is a bounded, linear map satisfying $X(p \cdot f)=p \cdot (Xf)$ for all $p  \in \mathbb{C}[z]$ and $f \in \clh$. Two Hilbert modules are said to be isomorphic if there exists a unitary module map between them.
\end{Definition}

Since one can extend the module action of a contractive Hilbert module $\clh$ over $\mathbb{C}[z]$ from
$\mathbb{C}[z]$ to the disk algebra $A(\mathbb{D})$ using the von
Neumann inequality, a contraction operator gives rise to a
contractive Hilbert module over $A(\mathbb{D})$. Recall that
$A(\mathbb{D})$ denotes the disk algebra, the algebra of
holomorphic functions on $\mathbb{D}$ that are continuous on the
closure of $\mathbb{D}$. Thus, the unitary equivalence of
contraction operators is the same as the isomorphism of the
associated contractive Hilbert modules over $A(\mathbb{D})$.

Next, let us recall that the Hardy space $H^2$ consists of
the holomorphic functions $f$ on $\mathbb{D}$ such that
$$\|f\|^2_2 = \sup_{0 <r <1} \frac{1}{2\pi} \int^{2\pi}_0 |f(re^{i\theta})|^2d\theta < \infty.$$ Similarly, the weighted Bergman
spaces $A^2_{\alpha}$, $-1 <
\alpha < \infty$, consist of the holomorphic functions $f$ on
$\mathbb{D}$ for which
$$\|f\|^2_{2, \alpha} = \frac{1}{ \pi} \int_{\mathbb{D}} |f(z)|^2
dA_{\alpha}(z) < \infty,$$ where $dA_{\alpha}(z) =
(\alpha +1) (1 - |z|^2)^{\alpha} dA(z)$ and $dA(z)$ denote the
weighted area measure and the area measure on $\mathbb{D}$,
respectively. Note that $\alpha=0$ gives the $($unweighted$)$ Bergman
space $A^2$. We mention \cite{Z} for a comprehensive treatment of the theory of Bergman spaces. The Hardy space, the Bergman space and the weighted Bergman spaces are contractive modules under the multiplication by the coordinate function.

The Hardy, the Bergman, and the weighted Bergman modules serve as examples of {\it
contractive reproducing kernel Hilbert modules}. A reproducing kernel Hilbert module
is a Hilbert module with a function called a {\it positive definite kernel} whose definition we now review.

\begin{Definition}
We say that a function $K : \mathbb{D} \times \mathbb{D} \raro \cll(\cle)$ for
a Hilbert space $\cle$, is a positive definite
kernel if
$\langle \sum_{i, j = 1}^{p} K(z_i, z_j)
\eta_i, \eta_j \rangle \geq 0$ for all $z_i
\in \mathbb{D}, \, \eta_i \in \cle$, and $p \in \mathbb{N}$.

\end{Definition}

Given a positive definite kernel $K$, we can construct a Hilbert
space $\clh_K$ of $\cle$-valued functions defined to be
$$\vee_{z \in \mathbb{D}} \vee_{ \eta \in \cle} K(\cdot, z) \eta ,$$ with inner product $$\langle K(\cdot, w) \eta,
\,K(\cdot, z) \zeta\rangle_{\clh_K} = \langle K(z, w) \eta,
\zeta\rangle_{\cle},$$ for all $z, w \in \mathbb{D}$ and $\eta,
\zeta \in \cle$. The evaluation of $f \in \clh_K$ at a point $z
\in \mathbb{D}$ is given by the reproducing property so that
$$\langle f(z), \eta \rangle_{\cle} = \langle f, K(\cdot, z)
\eta\rangle_{\clh_K},$$ for all $f \in \clh_K, z  \in \mathbb{D}$
and $\eta \in \cle$. In particular, the evaluation operator
$\bm{ev}_z : \clh_K \raro \cle$, $\bm{ev}_z (f) := f(z)$
is bounded for all $z \in \mathbb{D}$.

Conversely, given a Hilbert space $\clh$ of holomorphic
$\cle$-valued functions on $\mathbb{D}$ with bounded evaluation
operator $\bm{ev}_z \in \cll(\clh, \cle)$ for each $z \in
\mathbb{D}$, we can construct a reproducing kernel $$ \bm{ev}_z
\circ \bm{ev}_w^* : \mathbb{D} \times \mathbb{D} \raro
\cll(\cle),$$ for all $z, w \in \mathbb{D}$. To ensure that
$\bm{ev}_z \circ \bm{ev}_w^*$ is injective, we must assume for
every $z \in \mathbb{D}$ that $\overline{\{f(z): f \in \clh
\}}=\cle$.

 A reproducing kernel Hilbert module is said to be a {\it contractive reproducing kernel Hilbert module} over $A(\mathbb{D})$
 if the operator $M_z$ is contractive.

The kernel function for $H^2$ is $K(z, w) = (1 - \bar{w} z)^{-1}$.
For $A^2_{\alpha}$, it is $$K(z, w) = (1 - \bar{w} z)^{-2 -
\alpha} = \sum_{k=0}^{\infty} \frac{\Gamma(k+2+\alpha)}{k!
\Gamma(2+ \alpha)}(\bar{w} z)^k,$$ where $\Gamma$ is the gamma
function.

It is well known that the multiplier algebra of $\clh$ is $H^{\infty}$, that is, $M_{\varphi} \clh \subseteq \clh$, for $M_{\vp}$ the operator of multiplication by $\varphi \in H^{\infty}$, where $H^{\infty}=H^{\infty}(\mathbb{D})$ is the algebra of bounded, analytic functions on $\mathbb{D}$ and $\clh$ is $H^2$, $A^2$ or $A^2_{\alpha}$. Moreover, for all $w \in \mathbb{D}$ and $\varphi \in
H^{\infty}$, $$M^*_{\varphi} K(\cdot, w) =
 \overline{\varphi(w)} K(\cdot, w).$$

\newsection{The class $B_n(\mathbb{D})$}

In \cite{CD}, M. Cowen and the first author introduced a class of
operators $B_n(\mathbb{D})$, which includes $M_z^*$ for the
operator $M_z$ defined on contractive reproducing kernel Hilbert
modules of interest in this note. We now recall the notion of
$B_n(\mathbb{D})$. Let $\clh$ be a Hilbert space and $n$ a positive integer.

\begin{Definition}\label{20}
An operator $T \in \cll(\clh)$ is in
the class $B_n(\mathbb{D})$ if

(i) $\mbox{dim} \, \mbox{ker } (T - w) = n$ for all $w \in \mathbb{D}$,

(ii) $\vee_{w \in \mathbb{D}} \mbox{ker}\, (T - w) = \clh$, and

(iii) $\mbox{ran}\, (T - w)  = \clh$ for all $w \in \mathbb{D}$.

\end{Definition}

\begin{Remark}
Since it follows from (iii) that $T-w$ is semi-Fredholm for all $w
\in \mathbb{D}$, (iii) actually implies (i) if we assume that $\mbox{dim\, ker} (T - w) < \infty$
for some $w \in \mathbb{D}$.
\end{Remark}

It is a result of Shubin \cite{S} that for $T \in
B_n(\mathbb{D})$, there exists a hermitian holomorphic rank $n$
vector bundle $E_T$ over $\mathbb{D}$ defined as the pull-back of
the holomorphic map $w \mapsto \mbox{ker}\, (T - w)$ from
$\mathbb{D}$ to the Grassmannian $Gr(n, \clh)$ of the
$n$-dimensional subspaces of $\clh$. As mentioned earlier in the
Introduction, in this note we consider contraction operators $T$
such that $T^* \in B_n(\mathbb{D})$. In other words, we
investigate contractive Hilbert modules $\clh$ with $M_z^* \in
B_n(\mathbb{D})$. For simplicity of notation, we will write $\clh
\in B_n(\mathbb{D})$. Thus, we have an anti-holomorphic map $w
\mapsto \mbox{ker}\, (M_z - w)^*$ instead of a holomorphic one and
therefore obtain a frame $\{\psi_i\}_{i=1}^n$ of anti-holomorphic
$\clh$-valued functions on $\mathbb{D}$ such that
$$\vee_{i=1}^n \psi_i(w) = \mbox{ker}\, (M_z - w)^* \subseteq \clh,$$ for every $w \in \mathbb{D}$. We will use the
notation $E_{\clh}^*$ for this anti-holomorphic vector bundle since it is the dual of the natural
hermitian holomorphic vector bundle $E_{\clh}$ defined by localization.

One can show for an operator belonging to a ``weaker" class than $B_n(\mathbb{D})$ that there still exists an anti-holomorphic frame. Since having such a frame is sufficient for many purposes, one can consider operators in this ``weaker" class, which will be introduced after the following proposition:

\begin{Proposition}\label{frame}
  Let $T \in \cll(\clh)$ and $\tilde{T} \in \cll(\tilde{\clh})$. Suppose that there exist anti-holomorphic functions $\{\psi_i\}_{i=1}^n$ and $\{\tilde{\psi_i}\}_{i=1}^n$ from $\mathbb{D}$ to $\clh$ and $\tilde{\clh}$, respectively, satisfying

(1) $T \psi_i(w) = \bar{w} \psi_i(w)$ and $\tilde{T} \tilde{\psi}_i(w) = \bar{w} \tilde{\psi}_i(w),$ for all $1 \leq i \leq n$, $w \in \mathbb{D}$, and

(2) $\vee_{w \in \mathbb{D}} \vee_{i=1}^n \psi_i(w) = \clh$ and $\vee_{w \in \mathbb{D}} \vee_{i=1}^n \tilde{\psi}_i(w) = \tilde{\clh}$.

\NI Then there is an anti-holomorphic partial isometry-valued function $V(w) : \clh \raro \tilde{\clh}$ such that $\mbox{ker}\, V(w) = [ \vee_{i=1}^n \psi_i(w)]^{\perp}$ and $\mbox{ran}\, V(w) = \vee _{i=1}^n \tilde{\psi}_i(w) $ if and only if there exists a unitary operator $V : \clh \raro \tilde{\clh}$ such that $(V \psi_i) (w) = V(w) \psi_i(w)$ for every $1 \leq i \leq n$ and $w \in \mathbb{D}$.
\end{Proposition}

\NI \textsf{Proof.} We refer the reader to the proof of the rigidity theorem in \cite{CD}, where the language of bundles is used. \qed

It was pointed out by N. K. Nikolski to the first author that the basic calculation used to prove the rigidity theorem \cite{CD} appeared earlier in \cite{P}.

\begin{Definition}
Suppose $T \in \cll(\clh)$ is such that $\mbox{dim}\, \mbox{ker}\, (T - w) \geqslant n$ for all $w \in \mathbb{D}$. We say that $T$
is in the class $B_n^w(\mathbb{D})$ or weak-$B_n(\mathbb{D})$ if there exist anti-holomorphic functions $\{\psi_i\}_{i=1}^n$ from $\mathbb{D}$ to $\clh$ such that

(i) $\{\psi_i(w)\}_{i=1}^n$ is linearly independent for all $w \in \mathbb{D}$,

(ii) $\vee_{i=1}^n \psi_i(w) \subseteq \mbox{ker}\, (T -w)$ for all $w \in \mathbb{D}$, and

(iii) $\vee_{w \in \mathbb{D}} \vee_{i=1}^n \psi_i (w)  = \clh$.

\end{Definition}

\begin{Remark}
The class $B_n^w(\mathbb{D})$ is closely related to the one considered by Uchiyama in \cite{U}.
\end{Remark}

Since the $\{\psi_i\}_{i=1}^n$ in Definition 3.4 frame a rank $n$
hermitian anti-holomorphic bundle, it suffices for our purpose to consider
contractive Hilbert modules $\clh$ with $M_z^* \in
B_n^w(\mathbb{D})$ instead of those with $M^*_z \in
B_n(\mathbb{D})$. We will write $\clh \in B_n^w(\mathbb{D})$ to
represent this case.

We continue this section with a brief discussion of some complex geometric notions. Since the anti-holomorphic vector bundle  $E_{\clh}^*$ also has hermitian structure, one can define the canonical Chern connection $\cld_{E^*_{\clh}}$ on $E_{\clh}^*$ along with its associated curvature two-form $\clk_{E^*_{\clh}}$. For the case $n =1$, $E_{\clh}^*$ is a line bundle and

\begin{equation} \label{curvE} \mathcal K_{E^*_{\clh}}(z)  = - \frac{1}{4} \bigtriangledown^2 \log \|\gamma_z\|^2 \,dz \wedge d\bar{z},\end{equation} for $z \in \mathbb{D}$, where $\gamma_z$ is an anti-holomorphic cross section of the bundle. For instance, by taking $\gamma_z$ to be the kernel functions for $H^2$ and $A^2_{\alpha}$, we see that
$$\clk_{E^*_{H^2}}(z)
= - \frac{1}{(1 - |z|^2)^{2}},$$ and
$$
\clk_{E^*_{A^2_{\alpha}}}(z)
= - \frac{2+\alpha}{(1 - |z|^2)^{2}}.
$$

In \cite{CD}, M. Cowen and the first author proved that the curvature is a complete unitary invariant, that is, two Hilbert modules $\clh$ and $\tilde{\clh}$ in $B_1(\mathbb{D})$ are isomorphic if and only if for every $z \in \mathbb{D}$, $$\clk_{E^*_{\clh}}(z) = \clk_{E^*_{\tilde{\clh}}}(z).$$ Now that we have Proposition \ref{frame} available, the result can be extended to Hilbert modules in $B_1^w(\mathbb{D})$. Note that two weighted Bergman modules cannot be isomorphic to each another, that is, $A^2_{\alpha}$ is isomorphic to $A^2_{\beta}$ if and only if $\alpha = \beta$. We also conclude that the Hardy module $H^2$ cannot be isomorphic to the weighted Bergman modules $A^2_{\alpha}$.

\newsection{Proof of the main results}

Let $\Theta= \{ \theta_1, \theta_2 \} \in H^{\infty}_{\cll(\mathbb{C}, \mathbb{C}^2)}$ satisfy the corona condition. Now denote by $\clh_{\Theta}$ the quotient Hilbert module $(\clh \otimes \mathbb{C}^2) / \Theta \clh$, where $\clh$ is $H^2$, $A^2$, or $A^2_{\alpha}$. This means that we have the following short exact sequence
$$0 \longrightarrow \clh  \otimes \mathbb{C} \stackrel{M_{\Theta}}\longrightarrow
\clh \otimes \mathbb{C}^2 \stackrel{\pi_{\Theta}}\longrightarrow
\clh_{\Theta} \longrightarrow 0,$$ where the first map
$M_{\Theta}$ is $M_{\Theta} f = \theta_1 f \otimes e_1 +\theta_2 f
\otimes e_2$ and the second map $\pi_{\Theta}$ is the quotient
Hilbert module map. The fact that $\Theta$ satisfies the corona
condition implies that $\mbox{ran } M_{\Theta}$ is closed. We denote the module multiplication $P_{\clh_{\Theta}} (M_z \otimes I_{\mathbb{C}^2})|_{\clh_{\Theta}}$ of the quotient Hilbert module $\clh_{\Theta}$ by $N_z$. We will see later that $\clh_{\Theta} \in
B_1(\mathbb{D})$, but for the time being, we first show that
$\clh_{\Theta} \in B_1^w(\mathbb{D})$.

\begin{Theorem}\label{hardy-section}
For $\Theta=\{\theta_1, \theta_2 \}$ satisfying the corona condition, $\clh_{\Theta} \in B_1^w(\mathbb{D})$.

\end{Theorem}
\NI \textsf{Proof.}
We first prove that $\mbox{dim}\,\mbox{ker}\, (N_z-w)^* = 1$ for all $w \in \mathbb{D}$. To this end, let $I_w := \{ p(z) \in \mathbb{C}[z] : p(w) = 0\}$, a maximal ideal in $\mathbb{C}[z]$. One considers the localization of the sequence
$$0 \raro \clh \otimes \mathbb{C} \stackrel{M_\Theta} \raro  \clh \otimes \mathbb{C}^2 \stackrel{\pi_{\Theta}} \raro
\clh_{\Theta} \raro 0,$$ to $w \in \mathbb{D}$ to obtain
$$\clh / I_w \cdot \clh \longrightarrow (\clh \otimes \mathbb{C}^2) / I_w \cdot (\clh \otimes \mathbb{C}^2) \longrightarrow \clh_{\Theta}/ I_w \cdot \clh_{\Theta} \longrightarrow 0,$$
or equivalently, $$\mathbb{C}_w \otimes \mathbb{C} \stackrel{I_{\mathbb{C}_w} \otimes \Theta(w)} \longrightarrow \mathbb{C}_w \otimes \mathbb{C}^2 \stackrel{\pi_{\Theta}(w)} \longrightarrow \clh_{\Theta}/ I_w \cdot \clh_{\Theta} \longrightarrow 0.$$ Since this sequence is exact and $\mbox{dim } \mbox{ran } \Theta(w)=1$ for all $w \in \mathbb{D}$, we have $\mbox{dim } \mbox{ker } \pi_{\Theta}(w) = 1$ (see \cite{DP}). Thus, $\mbox{dim}\, \clh_{\Theta}/ I_w \cdot \clh_{\Theta} = 1$, and so $\mbox{dim }\mbox{ker }(N_z-w)^*=1$ for all $w \in \mathbb{D}$.

Now denote by $k_w$ a kernel function $k(\cdot, w)$ for $\clh$, and by $\{e_1, e_2\}$ an orthonormal basis for $\mathbb{C}^2$. We prove that $$\gamma_w := k_w \otimes (\overline{\theta_2(w)} e_1 - \overline{\theta_1(w)} e_2)$$ is a non-vanishing anti-holomorphic function from $\mathbb{D}$ to $\clh \otimes \mathbb{C}^2$ such that

(1) $\gamma_w \in \mbox{ker}\, (N_z-w)^*$ for all $w \in \mathbb{D}$, and

(2) $\vee_{w \in \mathbb{D}} \gamma_w = \clh_{\Theta}$.

Since the $\theta_i$ are holomorphic and $k_w$ is anti-holomorphic, the fact that $w \mapsto \gamma_w$ is anti-holomorphic follows. Furthermore, since $\Theta$ satisfies the corona condition, the $\theta_i$ have no common zero and hence $\gamma_w  \neq \bm{0}$ for all $w \in \mathbb{D}$. Now, for $f \in  \clh$, $M_{\Theta} f = \theta_1 f \otimes e_1 + \theta_2 f \otimes e_2$ and therefore for all $w \in \mathbb{D}$,

\[
\begin{split}
\langle M_{\Theta} f, \gamma_w \rangle & = \langle \theta_1 f, k_w \rangle \langle e_1, \overline{\theta_2(w)} e_1 \rangle - \langle \theta_2 f, k_w \rangle \langle e_2, \overline{\theta_1(w)} e_2 \rangle \\& = \theta_1(w)f(w)\theta_2(w)-\theta_2(w)f(w)\theta_1(w)=0.
\end{split}
\]

Hence, $\gamma_w \in (\mbox{ran} M_{\Theta})^{\perp} = \clh_{\Theta}$. Moreover, since $M_z^*k_w=\bar{w}k_w$ for all $w \in \mathbb{D}$,
\[
\begin{split}
N_z^* \gamma_w & = (M_z \otimes I_{\mathbb{C}^2})^* \gamma_w = M_z^* (\overline{\theta_2(w)} k_w) \otimes e_1 - M^*_z (\overline{\theta_1(w)} k_w) \otimes e_2 \\& = \overline{\theta_2(w)} \bar{w} k_w \otimes e_1 - \overline{\theta_1(w)} \bar{w} k_w \otimes e_2\\& = \bar{w} \gamma_w.
\end{split}
\]

Next, in order to show that (2) holds, it suffices to prove that
for $h = h_1 \otimes e_1 + h_2 \otimes e_2 \in \clh \otimes
\mathbb{C}^2$ such that $h \perp \vee_{w \in \mathbb{D}}
\gamma_w$, we have $h \in \mbox{ran}\ M_\Theta$. We first claim
that there exists a function $\eta$ defined on
$\mathbb{D}$ such that for all $w \in \mathbb{D}$ and $i=1,2$,
$$h_i(w) =   \theta_{i}(w) \eta(w).$$
Since $h \perp \gamma_w$ for  every $w \in \mathbb{D}$, we have
\[
\begin{split}
\langle h, \gamma_w \rangle & = \langle h_1, k_w \rangle \langle
e_1, \overline{\theta_2(w)} e_1 \rangle - \langle h_2, k_w \rangle
\langle e_2, \overline{\theta_1(w)} e_2 \rangle \\& =
h_1(w)\theta_2(w)-h_2(w)\theta_1(w)=0,
\end{split}
\]
or equivalently,
\begin{equation}
\mbox{det~}\ \begin{bmatrix} h_1(w) & \theta_1(w) \\ h_2(w) &
\theta_2(w)
\end{bmatrix}=0,
\end{equation}
for all $w \in \mathbb{D}$. Thus using the fact that $\mbox{rank
}\begin{bmatrix} \theta_1(w) \\ \theta_2(w)
\end{bmatrix}=1$ for all $w \in \mathbb{D}$, we obtain a unique nonzero function $\eta(w)$
satisfying $h_i(w)=\theta_i(w) \eta(w)$ for $i=1,2$.

The proof is completed once we show that $\eta \in \clh$. Note
that by the corona theorem, we get $\psi_1, \psi_2 \in
H^{\infty}$ such that $\psi_1(w) \theta_1(w)+\psi_2(w)
\theta_2(w)=1$ for every $w \in \mathbb{D}$. Since $\eta=(\psi_1
\theta_1 + \psi_2 \theta_2)\eta=\psi_1 h_1 + \psi_2 h_2$, and
$H^{\infty}$ is the multiplier algebra for $\clh$, the result
follows. \qed

\begin{Remark}
\textnormal{Observe that the above proof shows that the hermitian
anti-holomorphic line bundle corresponding to the quotient Hilbert
module $\clh_{\Theta}$ is the twisted vector bundle obtained as
the bundle tensor product of the hermitian anti-holomorphic line
bundle for $\clh$ with the anti-holomorphic dual of the line
bundle $\coprod_{w \in \mathbb{D}} \mathbb{C}^2/ \Theta(w)
\mathbb{C}$. This phenomenon holds in general; suppose that for
Hilbert spaces $\cle$ and $\cle_*$, $\Theta \in
H^{\infty}_{\cll(\cle, \cle_*)}$ and $M_{\Theta}$ has closed
range. If the quotient Hilbert module $\clh_{\Theta}$,
$$0  \raro \clh \otimes \cle \stackrel{M_\Theta} \raro \clh \otimes \cle_* \raro \clh_{\Theta} \raro 0,$$ is in $B_n(\mathbb{D})$,
then the rank $n$ hermitian anti-holomorphic vector bundle
$E^*_{\clh_{\Theta}}$ for $\clh_{\Theta}$ is the bundle tensor
product of $E^*_{\clh}$ with the anti-holomorphic dual of the rank
$n$ bundle $\coprod_{w \in \mathbb{D}} \cle_*/{\Theta(w) \cle}$ (see \cite{DKKS}).}
\end{Remark}

In order to have $\clh_{\Theta} \in B_1(\mathbb{D})$, it now remains to check only one condition. We do this in the following Proposition.

\begin{Proposition}
$\mbox{ran }(N_z-w)^*=\clh_\Theta$ for all $w \in \mathbb{D}$.
\end{Proposition}
\NI \textsf{Proof.} We write
\begin{equation*}
M_z \otimes I_{\mathbb{C}^2}\backsim \begin{bmatrix} * & *\\0 & N_z\end{bmatrix}
\end{equation*}
relative to the decomposition $\clh \otimes \mathbb{C}^2 = \mbox{ran } M_{\Theta} \oplus (\mbox{ran}\, M_{\Theta})^{\perp}$. It suffices to note that $\clh \in B_1(\mathbb{D})$ implies that $\mbox{ran }(M_z-w)^*=\clh$. \qed

Let us now consider the curvature $\clk_{E^*_{{\clh}_{\Theta}}}$.
By (3.1), one needs only to compute the norm of the section
$$\gamma_w = k_w \otimes
(\overline{\theta_2(w)}e_1-\overline{\theta_1(w)}e_2)$$ given in
Theorem 4.1. Since $$\|\gamma_w\|^2 = \|k_w\|^2 (|\theta_1(w)|^2 +
|\theta_2(w)|^2),$$ we get the identity
\begin{equation}
\clk_{E^*_{{\clh}_{\Theta}}}(w)=\clk
_{E^*_{\clh}}(w)-\frac{1}{4} \bigtriangledown^2
\mbox{log}\,(|\theta_1(w)|^2 + |\theta_2(w)|^2),
\end{equation}
for all $w \in \mathbb{D}$.

We are now ready to prove Theorem 4.4. For the sake of convenience, we restate it here.

\begin{Theorem}
Let $\Theta=\{\theta_1, \theta_2\}$ and $\Phi=\{\vp_1, \vp_2\}$ satisfy
the corona condition. The quotient Hilbert
modules ${\clh}_{\Theta}$ and ${\clh}_{\Phi}$ are isomorphic if and only
if
$$
\bigtriangledown^2 \mbox{log}\,\frac{|\theta_1(z)|^2 +
|\theta_2(z)|^2}{|\vp_1(z)|^2 + |\vp_2(z)|^2}=0,
$$
for all $z \in \mathbb{D}$, where $\clh$ is the Hardy, the Bergman, or a weighted Bergman module.
\end{Theorem}

 \NI \textsf{Proof.} Since ${\clh}_{\Theta}, {\clh}_{\Phi} \in B_1^w(\mathbb{D})$, (we have seen that they actually belong to $B_1(\mathbb{D})$), they are isomorphic if and only if $\clk_{E^*_{{\clh}_{\Theta}}}(w)=\clk_{E^*_{{\clh}_{\Phi}}}(w)$ for all $w \in \mathbb{D}$. But note that (4.2) and an analogous identity for $\Phi$ hold, where the $\theta_i$ are replaced with the $\vp_i$. Since both $\Theta$ and $\Phi$ satisfy the corona condition, the result then follows.
\qed \\

We once again state Theorem 4.5.

\begin{Theorem}
Suppose that $\Theta=\{\theta_1, \theta_2\}$ and $\Phi=\{\vp_1, \vp_2\}$ satisfy
the corona condition. The quotient Hilbert modules
$(A^2_{\alpha})_{\Theta}$ and
$(A^2_{\beta})_{\Phi}$ are isomorphic if and
only if $\alpha=\beta$ and
\begin{equation}
\bigtriangledown^2 \mbox{log}\,\frac{|\theta_1(z)|^2 +
|\theta_2(z)|^2}{|\vp_1(z)|^2 + |\vp_2(z)|^2}=0,
\end{equation}
for all $z \in \mathbb{D}$.
\end{Theorem}

 \NI \textsf{Proof.} Since we have by (4.2), $$\clk_{E^*_{(A^2_\alpha)_\Theta}}(w) =  - \frac{2+\alpha}{(1 - |w|^2)^{2}} - \frac{1}{4} \bigtriangledown^2 \mbox{log}\,(|\theta_1(w)|^2 + |\theta_2(w)|^2),$$ and $$\clk_{E^*_{(A^2_\beta)_\Phi}} (w) =  - \frac{2+\beta}{(1 - |w|^2)^{2}} - \frac{1}{4} \bigtriangledown^2 \mbox{log}\,(|\vp_1(w)|^2 + |\vp_2(w)|^2),$$ one implication is obvious.
 For the other one, suppose that $(A^2_ \alpha)_{\Theta}$ is isomorphic to $(A^2_\beta)_{\Phi}$ so that
 the curvatures coincide. Observe next that $$\frac{4(\beta - \alpha)}{(1 - |w|^2)^{2}} =  \bigtriangledown^2 \mbox{log}\,\frac{|\theta_1(w)|^2 + |\theta_2(w)|^2}{|\vp_1(w)|^2 + |\vp_2(w)|^2}.$$
 Since a function $f$ with $\bigtriangledown^2 f(z) = \frac{1}{(1 - |z|^2)^2}$ for all $z \in \mathbb{D}$ is necessarily unbounded, we have a contradiction unless $\alpha=\beta$ (see Lemma 4.6 below) and (4.3) holds. This is due to the assumption that the bounded functions $\Theta$ and $\Phi$ satisfy the corona condition. \qed

 \begin{Lemma} There is no bounded function $f$ defined on the unit disk $\mathbb{D}$ that satisfies $\bigtriangledown^2 f(z) =\frac{1}{(1 - |z|^2)^2}$ for all $z \in \mathbb{D}$.
 \end{Lemma}

\NI \textsf{Proof. } Suppose that such $f$ exists. Since $\frac{1}{4} \bigtriangledown^2 [(|z|^2)^{m}] = \partial \bar{\partial} [(|z|^2)^{m}] = m^2 (|z|^2)^{m-1}$ for all $m \in \mathbb{N}$, we see that for $$g(z) := \frac{1}{4} \sum_{m=1}^{\infty} \frac{|z|^{2m}}{m} = - \frac{1}{4} \mbox{log}\, (1 - |z|^2),$$ $ \bigtriangledown^2 g(z) = \frac{1}{(1-|z|^2)^2}$ for all $z \in \mathbb{D}$. Consequently, $f(z)=g(z)+h(z)$ for some harmonic function $h$. Since the assumption is that $f$ is bounded, there exists an $M>0$ such that $|g(z) + h(z)| \leq M$ for all $z \in \mathbb{D}$. It follows that $$ \mbox{exp}\,(h(z))  \leq \mbox{exp}\,( - g(z) + M) = (1 - |z|^2)^{\frac{1}{4}} \,\mbox{exp}\,(M),$$ and letting $z = r e^{i\theta},$ we have $\mbox{exp}\,(h(r e^{i\theta})) \leq (1 - r^2)^{\frac{1}{4}} \mbox{exp}\,(M)$. Thus $\mbox{exp}\,(h(r e^{i\theta})) \raro 0$ uniformly as $r \raro 1^-$, and hence $\exp h(z) \equiv 0$. This is due to the maximum modulus principle because $\exp h(z) = | \exp(h(z) + i \tilde{h}(z))|$, where $\tilde{h}$ is a harmonic conjugate for $h$. We then have a contradiction, and the proof is complete. \qed \\

We thank E. Straube for providing us with a key idea used in the proof of Lemma 4.6.

\begin{Theorem}
For $\Theta=\{\theta_1, \theta_2\}$ and $\Phi=\{\vp_1, \vp_2\}$ satisfying
the corona condition, $(H^2)_{\Theta}$ cannot be isomorphic to $(A^2_{\alpha})_{\Phi}$.
\end{Theorem}

\NI \textsf{Proof.} By identity (4.2), we conclude that $(H^2)_{\Theta}$ is isomorphic to $(A^2_{\alpha})_{\Phi}$ if and only if $$\frac{4(1+\alpha)}{(1 - |w|^2)^{2}} =  \bigtriangledown^2 \mbox{log}\,\frac{|\vp_1(w)|^2 + |\vp_2(w)|^2}{|\theta_1(w)|^2 + |\theta_2(w)|^2}.$$ But according to Lemma 4.6, this is impossible unless $\alpha=-1$.\qed

\newsection{Concluding remark}
Although the case of quotient modules we have been studying in this note may seem rather elementary, the class of examples obtained is not without interest. The ability to control the data in the construction, that is, the multiplier, provides one with the possibility of obtaining examples of Hilbert modules over $\mathbb{C}[z]$ and hence operators with precise and refined properties. In \cite{BDF} and \cite{BDFP} the authors utilized this framework to exhibit operators with properties that responded to questions raised in the papers.

In particular, in \cite{BDFP} the authors are interested in characterizing contraction operators that are quasi-similar to the unilateral shift of multiplicity one. In the earlier part of the paper, which explores a new class of operators, a plausible conjecture presents itself but examples defined in the framework of this note, introduced in Corollary 7.9, show that it is false.

In \cite{BDF}, the authors study canonical models for bi-shifts; that is, for commuting pairs of pure isometries. A question arises concerning the possible structure of such pairs and again, examples built using the framework of this note answer the question.

Finally in \cite{DMS}, the authors determine when a contractive Hilbert module in $B_1(\mathbb{D})$ can be represented as a quotient Hilbert module of the form $\clh_\Theta$, where $\clh$ is the Hardy, the Bergman, or a weighted Bergman module. For the case of the Hardy module, the result is contained in the model theory of Sz.-Nagy and Foias \cite{NF}.

One can consider a much larger class of quotient Hilbert modules replacing the Hardy, the Bergman and the weighted Bergman modules by a quasi-free Hilbert module \cite{DM1} of rank one. In that situation, one can raise several questions relating curvature invariant, similarity and the multiplier corresponding to the given quotient Hilbert modules. These issues will be discussed in the forthcoming paper \cite{DKKS}.

\end{document}